\documentclass[12pt]{article} 
\begin{document}

\title{Uniform partitions of 3-space, their relatives and embedding
\thanks{This work was supported by the Volkswagen-Stiftung 
(RiP-program at Oberwolfach) and Russian fund of fundamental 
research (grant 98-01-00251).}}

\author{Michel DEZA \\
     CNRS and Ecole Normale Superieure, Paris, France \\
\and   Mikhail SHTOGRIN \\
Steklov Mathematical Institute, 117966 Moscow GSP-1, Russia} 

\date{}

\maketitle 

\begin{abstract}
We review 28 uniform partitions of 3-space in order to find out which
of them have graphs (skeletons)  embeddable isometrically (or with
scale 2) into some cubic lattice ${\bf Z}_n$. We also consider some
relatives of those 28 partitions, including Achimedean 4-polytopes of
Conway-Guy,  non-compact uniform partitions, Kelvin partitions and
those with unique vertex figure (i.e. Delaunay star). Among last ones we
indicate two continuums of aperiodic tilings by semi-regular 3-prisms
with cubes or with regular tetrahedra and regular octahedra.
On the way many new partitions are added to incomplete cases
considered here.

\end{abstract}
\section{Introduction}
A polyhedron is called {\em uniform} if all its faces are regular polygons
and its group of symmetry is vertex-transitive.
A normal partition of 3-space is called {\em uniform} if all its facets (cells)
are uniform polyhedra and group of symmetry is vertex-transitive.
There are exactly 28 uniform partitions of 3-space.
A short history of this result follows. Andreini in 1905 proposed, as the
complete list, 25 such partitions. But one of them ($13'$, in his notation)
turns out to be not uniform; it seems, that Coxeter \cite{Cox}, page 334 was
the first to realize it. Also Andreini missed partitions 25-28 (in our
numeration given below). Till recent years, mathematical literature was
abundant with incomplete lists of those partitions. See, for example, 
\cite{Cr}, \cite{Wi} and \cite{Pe} (all of them does not contain
24-28) and \cite{Ga}. 
The first to publish the complete list was Gr\"unbaum in \cite{Gr}. 
But he wrote there that, after obtaining the list, he realized that 
the manuscript \cite{Jo} already contained all 28 partitions.
We also obtained all 28 partitions independently, but only in 1996.

We say that given partition $P$ has an {\em $l_1$-graph} and embeds up to
{\em scale} $\lambda$ 
into the cubic lattice ${\bf Z}_m$, if there exists a mapping $f$ of 
the vertex-set of the skeleton of $P$ into ${\bf Z}_m$ such that
\[\lambda d_{P}(v_i,v_j)=||f(v_i),f(v_j)||_{l_1} =
\sum_{1 \le k \le m}|f_{k}(v_i)-f_k(v_j)| \mbox{ for all vertices
$v_i,v_j$}. \]
The smallest such number $\lambda$ is called {\em minimal scale}; all
embeddings below are given with their minimal scale.

Call an $l_1$-partition $l_1$-{\em rigid}, if all its embeddings (as above) 
into cubic lattices are pairwise equivalent, i.e. unique up to a symmetry
of the cubic lattice.
All embeddable partitions (and all embeddable tiles, except of Tetrahedron, 
having two non-equivalent embeddings: into $\frac{1}{2}H_{3}$ and
$\frac{1}{2}H_{4}$) in this paper
turn out to be $l_1$-rigid and so, by a result from \cite{Shp}, having minimal
scale 1 or 2. Those embeddings were obtained by constructing a
complete system of {\em alternated zones}; see \cite{CDG}, \cite{DS1},
\cite{DS2}, \cite{DS3}.
 
The following {\em 5-gonal} inequality (\cite{De}):
 \[d_{xy}+(d_{ab}+d_{ac}+d_{bc}) \le
(d_{xa}+d_{xb}+d_{xc})+(d_{ya}+d_{yb}+d_{yc}) \]
 for distances between any five vertices $a,b,c,x,y$,
 is an important necessary condition for embedding of graphs.
It turns out that all non-embeddable partitions and tiles, considered in this
paper, are, moreover, not 5-gonal.

Denote by $De(T)$ and $Vo(T)$ the Delaunay and Voronoi partitions of 3-space
associated with given set of points $T$. By an abuse of language, we will
use the same notation for the graph, i.e. the skeleton of a partition.   
The Voronoi and Delaunay partitions are dual one to each other (not only
combinatorially, but metrically). Denote by $P^*$ the partition dual to
partition $P$; it should not be confounded with the same notation for
dual {\em lattice}.

In Tables 1, 2 we remind the results (from \cite{DS1}) on embedding
of uniform polyhedra and plane partitions.
In the Table \ref{tab:t1}, ${d(P)}$ denotes the diameter of polyhedron $P$.

\begin{table}
\caption{\bf Embedding of uniform polyhedra and their duals}
\label{tab:t1}
\begin{center}
\begin{tabular}{|c|c|c|c|c|} \hline
Polyhedron $P$; face-vector& $\ell_1$-status of $P$ &$d(P)$ &
 $\ell_1$-status of $P^{*}$ &
$d(P^*)$ \\ \hline \hline
{\em Tetrahedron; (3.3.3)}&$\to\frac{1}{2}H_3$&1&$\to\frac{1}{2}H_{4}$&1\\ \hline
{\em Cube; (4.4.4)}&$\to H_3$&3&$\to\frac{1}{2}H_{4}$&2\\ \hline
{\em Dodecahedron; (5.5.5)}&$\to\frac{1}{2}H_{10}$&5&$\to\frac{1}{2}H_{6}$&3\\ \hline \hline
{\em Cuboctahedron; (3.4.3.4)}  & non 5-gonal & 3 & $\to H_{4}$ & 4 \\ \hline
{\em Icosidodecahedron; (3.5.3.5)}  & non 5-gonal & 5 & $\to H_{6}$ & 6 \\ \hline
{\em Truncated tetrahedron; (3.6.6)}&non 5-gonal&3&$\to\frac{1}{2}H_{7}$&2\\ \hline
{\em Truncated octahedron; (4.6.6)}  & $\to H_{6}$ & 6 & non 5-gonal & 3\\ \hline
{\em Truncated cube; (3.8.8)} &non 5-gonal&6& $\to\frac{1}{2}H_{12}$& 3 \\ \hline
{\em Truncated icosahedron; (5.6.6)}&non 5-gonal&9&$\to\frac{1}{2}H_{10}$&5\\ \hline
{\em Truncated dodecahedron; (3.10.10)}&non 5-gonal&10& $\to\frac{1}{2}H_{26}$&5\\ \hline
{\em Rhombicuboctahedron; (3.4.4.4)}&$\to\frac{1}{2}H_{10}$&5&non 5-gonal&5\\ \hline
{\em Rhombicosidodecahedron; (3.4.5.4)}& $\to\frac{1}{2}H_{16}$ &8&non 5-gonal&8\\ \hline
{\em Truncated cuboctahedron; (4.6.8)}& $\to H_{9}$&9&non 5-gonal & 4\\ \hline
{\em Truncated icosidodecahedron; (4.6.10)}&$\to H_{15}$&15&non 5-gonal&6\\ \hline
{\em Snub cube; (3.3.3.3.4)}& $\to\frac{1}{2}H_{9}$&4&non 5-gonal&7 \\ \hline
{\em Snub dodecahedron; (3.3.3.3.5)}&$\to \frac{1}{2}H_{15}$&7&non 5-gonal&15\\ \hline
{\em $Prism_3$; (4.3.4)}&$\to\frac{1}{2}H_5$&2&$\to\frac{1}{2}H_{4}$&2\\ \hline
{\em $Prism_n$ (odd $n \ge 5$); (4.n.4)}&$\to\frac{1}{2}H_{n+2}$&
 $\lfloor \frac{n+2}{2} \rfloor$&non 5-gonal&2\\ \hline
{\em $Prism_n$ (even  $n \ge 6$); (4.n.4)}&$\to  H_{\frac{n+2}{2}}$&
 $\frac{n+2}{2}$&non 5-gonal&2\\ \hline
{\em $Antiprism_n$ ($n \ge 4$); (3.3.n.3)}&$\to\frac{1}{2}H_{n+1}$
 & $\lfloor \frac{n+1}{2} \rfloor$
 &non 5-gonal&3\\ \hline
\end{tabular}
\end{center}
\end{table}

\begin{table}
\caption{\bf Embedding of uniform plane tilings and their duals}
\label{tab:t2}
\begin{center}
\begin{tabular}
{|c|c|c|} \hline
Tiling $T$; face-vector& $\ell_1$-status of $T$ &
 $\ell_1$-status of $T^{*}$ \\ \hline \hline
{\em $Vo(Z_2)=De(Z_2)$; (4.4.4.4)}& $\to Z_{2}$& $\to Z_{2}$ \\ \hline
{\em $De(A_2)$; (3.3.3.3.3.3)}& $\to \frac{1}{2}Z_{3}$& $\to Z_{3}$ \\ \hline
{\em $Vo(A_2)$; (6.6.6)}& $\to Z_{3}$& $\to \frac{1}{2}Z_{3}$ \\ \hline \hline
{\em Kagome net; (3.6.3.6)}  & non 5-gonal & $\to Z_{3}$ \\ \hline
{\em (3.4.6.4)}  & $\to \frac{1}{2}Z_{3}$& non 5-gonal  \\ \hline
{\em truncated (4.4.4.4); (4.8.8)}&$\to Z_4$& non 5-gonal \\ \hline
{\em (4.6.12)}&$\to Z_6$& non 5-gonal \\ \hline
{\em truncated (6.6.6); (3.12.12)}&non 5-gonal&$\to \frac{1}{2}Z_{\infty}$
 \\ \hline
{\em chiral net; (3.3.3.3.6)}&$\to \frac{1}{2}Z_6$& non 5-gonal \\ \hline
{\em 3,3,3,4,4}&$\to \frac{1}{2}Z_3$& non 5-gonal \\ \hline
{\em dual Cairo net; (3.3.4.3.4)}&$\to \frac{1}{2}Z_4$& non 5-gonal \\ \hline

\end{tabular}
\end{center}
\end{table}

\newpage

\section{28 uniform partitions} 
In Table 3 of 28 partitions, the meaning of the column is:

1. the number which we give to the partition; 

2. its number in \cite{An} if any; 

3. its number in \cite{Gr}; 

4. a characterization (if any) of the partition;

5. tiles of partition and respective number of them in Delaunay star; 

5$^*$. tiles of its dual;

6. embeddability (if any) of partition; 

6$^*$. embeddability (if any) of its dual.

Notation $\frac{1}{2}{\bf Z}_m$ in columns 6, 6$^*$ means that
the embedding is isometric up to scale 2.

The notations for the tiles given in the Table 3 are:
tr$P$ for truncated polyhedron $P$; $Prism_n$ for semi-regular
n-prism; $\alpha_3$, $\beta_3$ and $\gamma_3$ for the Platonic 
tetrahedron, octahedron and cube; 
$Cbt$ and $Rcbt$ for Archimedean Cuboctahedron and
Rhombicuboctahedron; 
$RoDo$, tw$RoDo$ and $RoDo-v$ for Catalan Rhombic Dodecahedron, for
its twist and for $RoDo$ with deleted vertex of valency 3; $Pyr_4$ and
$BPyr_3$ for corresponding pyramid and bi-pyramid; $BDS^*$ for dual
bidisphenoid. 
   
Remark that \cite{Cox} considered 12 of all 28 partitions; namely,
No's 8, 7, 18, 2, 16, 23, 9 denoted there as $t_A \delta _4$ for 
$A= \{1\}$, $\{0,1\}$, $\{0,2\}$, $\{1,2\}$, $\{0,1,2\}$,
$\{0,1,3\}$, $\{0,1,2,3\}$, respectively, and No's 6, 5, 20, 19, 17 denoted as 
$q \delta_4$, $h \delta_4$, $h_{2} \delta_{4}$, $h_{3} \delta_4$, 
$h_{2,3}\delta_{4}$.

{\bf Table 3. Embedding of uniform partitions and their duals.}
\[\begin{array}{||l|l|l||l||lr|l||l|l||} \hline
1 & 2 & 3 & 4 & 5 && 5^* & 6 & 6^* \\ \hline
No & & &  &\mbox{tiles}& &\mbox{of dual}& \mbox{emb.}&\mbox{dual}\\
\hline        
1  &  1 &   22 &   De({\bf Z}_3)=Vo({\bf Z}_{3}) &\gamma_3 &8&\gamma_3 &
{\bf Z}_3  & {\bf Z}_3 \\
2  &  3 &  28  &  Vo( A_3^*=\mbox{bcc}) &tr\beta_3 &4& \sim \alpha_3 &
{\bf Z}_6 & - \\
3 &   4 &     11 &   De(A_2 \times {\bf Z}_1)  & Prism_3& 12&Prism_6 &
\frac{1}{2}{\bf Z}_4  & {\bf Z}_4 \\
4 &   5 &  26 &      Vo(A_2 \times {\bf Z}_1)  & Prism_6& 6&Prism_3 &
{\bf Z}_4  & \frac{1}{2}{\bf Z}_4 \\
5 &   2 &     1 &    De(A_3=\mbox{fcc}) & \alpha_3,\beta_3& 8,6& RoDo&  
\frac{1}{2}{\bf Z}_4 & {\bf Z}_4 \\ 
6 &   13 & 6 & \mbox{F\"oppl partition} &\alpha_3, tr \alpha_3& 2,6&\sim
\gamma_3 & -& {\bf Z}_4 \\
7 &  14  &   8 &  \mbox{boron CaB}_6  & \beta_3, tr \gamma_3& 1,4 
& Pyr_4 &  -  &  - \\
8 &   15 &    7 &    De(J\mbox{-complex})  &\beta_3, Cbt& 2,4& \sim 
\beta_3 & -  &     - \\
9 &   22 &    27 &        &Prism_8, & &\sim \alpha_3 & {\bf
Z}_9 &     - \\
 & & & & tr Cbt&2,2 & & & \\
10&   6 &     24 &   De((4.8^2 \times {\bf Z}_1) &Prism_8, \gamma_3& 4,2 
&\sim Prism_3 & {\bf Z}_5 &     - \\
11 &  8 &     18 &   De(3.6.3.6 \times {\bf Z}_1) & Prism_3, 
 & & \sim \gamma_3 &-  & {\bf Z}_4 \\ 
 & & & &Prism_6 & 4,4 & & & \\
12 &  12 &    17 &  De(3^4.6 \times {\bf Z}_1) & Prism_3, &  
&\sim Prism_5 & \frac{1}{2}{\bf Z}_7 &    - \\ 
 & & & &Prism_6 & 8,2 & & & \\
13 &  11 &    13 & De(3^3.4^2 \times {\bf Z}_1) & Prism_3, \gamma_3& 6,4
& \sim Prism_5 &\frac{1}{2}{\bf Z}_4 &    - \\ 
14 &  11'&    14 & De(3^2.4.3.4 \times {\bf Z}_1) & Prism_3,
\gamma_3& 6,4
& \sim Prism_5 &\frac{1}{2}{\bf Z}_5 &  - \\ 
15 &  7 &     19 & De(3.12^2 \times {\bf Z}_1) & Prism_3, & 
&\sim Prism_3 & -& \frac{1}{2}{\bf Z}_{\infty} \\
 & & & &Prism_{12} & 2,4 & & & \\
16 &  18 &    25 & \mbox{a zeolit}& \gamma_3, tr \beta_3, tr Cbt &
1,1,2& \sim \alpha_3 &{\bf Z}_9 &-\\
17 &  20 &    21 &       &tr \alpha_3, tr \gamma_3, & 
&\sim \alpha_3 &- &-\\
 & & & & tr Cbt & 1,1,2 & & & \\
18 &  17 &    9  &   &\gamma_3, Cbt, Rcbt& 2,12 &\sim BPyr_3 &- &-\\
19 &  16 &    5  &  & \alpha_3, \gamma_3, Rcbt& 1,3,1&\sim BPyr_3
&\frac{1}{2}{\bf Z}_7 & - \\ 
20 &  21 &    10 & \mbox{boron UB}_{20} &tr \alpha_3, tr \beta_3,
Cbt& 2,1,2
&\sim Pyr_4 &- &-\\
21 &  9  &    16 & De(3.4.6.4 \times {\bf Z}_1) & \gamma_3, 
Prism_3,& &\sim \gamma_3 &\frac{1}{2}{\bf Z}_4 &-\\ 
   &     &  &  & Prism_6 & 4,2,2 & & & \\
22 &  10 &    23 & De(4.6.12 \times {\bf Z}_1) &\gamma_3, 
Prism_3, & &\sim Prism_3 & {\bf Z}_7 &-\\
   &     &  &  & Prism_{12} &2,2,2  & & & \\
23 &  19 &    20 &   &\gamma_3, Prism_8, &  &\sim Pyr_4 &- &-\\
   &     &       &   & tr \gamma_3, Rcbt& 1,2,1,1 &           &  & \\  
24 &  2' &    2  & De(\mbox{hcp})  & \alpha_3, \beta_3& 8,6& tw RoDo & - & - \\
25 &  - &     3  & De(\mbox{elong. }A_3) &Prism_3, \beta_3, \alpha_3 &
6,4,3 & RoDo-v & \frac{1}{2}{\bf Z}_4 &{\bf Z}_4 \\
26 &  - &     4 & De(\mbox{elong. hcp}) &Prism_3, \beta_3,\alpha_3 &
6,4,3 & RoDo-v& -& - \\
27 &  - &     12 &   & Prism_3& 12& BDS^* & \frac{1}{2}{\bf Z}_5 & - \\
28 &  - &     15 & De(\mbox{elong.
 27}) & Prism_3, \gamma_3& 6,4 
&\sim Prism_5 & \frac{1}{2}{\bf Z}_5 & - \\ \hline
\end{array} \]

Remarks to Table 3: 

1. The partition 15$^*$ is only one embeddable into ${\bf Z}_{\infty}$ 
(in fact, with scale 2). 

2. All partitions embeddable with scale 1 are, except 25$^*$, zonohedral.
  The Voronoi tile of 25$^*$  is not centrally-symmetric. 
It will be interesting to find a normal tiling of 3-space embeddable
with scale 1 such that the tile is centrally-symmetric; such
non-normal tiling is given in 
\cite{Sh}: see item 35 in Table 4 below. 

3. An embedding of tiles of tiling is necessary, but not sufficient, for
embedding of whole tiling; for example, $26^*$ and $27^*$ are non-embeddable
while their tiles are embeddable into $H_4$ and
$\frac{1}{2}H_8$, respectively. In fact, all dual uniform partitions $P^*$
in Table 3, except of 24$^*$, have no non-embeddable tiles. Among all 11 non-embeddable uniform
partitions only items 11, 15, 24 and 26 have only embeddable tiles. The same
is true for tilings 30, 33 and 32$^*$, 33$^*$, 34$^*$, 46$^*$ of Table 4. 

4. Among all 28 partitions only No's 1, 2, 5, 6, 8  have same surrounding 
of edges: polygons (4.4.4.4), (4.6.6), (3.3.3.3), (3.3.6.6), (3.3.4). 

5. Partitions 8 and 24 are Delaunay partitions of lattice complexes: namely, a
3-lattice called J-complex and a bi-lattice hcp; the tile of $Vo$(J-complex)
has form of jackstone (it explains the term "J-complex")) and it is
combinatorially equivalent to $\beta_3$.

6. Partitions 1, 3, 5 are Delaunay partitions of lattices 
${\bf Z}_3$, $A_2 \times {\bf Z}_1$, $A_3$=fcc.
  Partitions 2 and 4 are Voronoi partitions of lattices $A_3^*$ =bcc 
  and $A_2 \times {\bf Z}_1$.
  No's 10, 11, 12, 13, 14, 15, 21, 22 are Delaunay prismatic partitions 
over 8 Archimedean partitions of the plane; the embeddability of them and their
duals is the same as in Table 2, but the dimension increases by 1.

7. Partitions 7 and 20 occur in Chemistry as borons CaB$_6$ and UB$_{20}$, 
respectively. Partition 20 occurs in zeolites.

8. The ratio of tiles in partition is 1:1 for 6, 7, 8, 10; 2:1 for 5,
11, 13, 14, 24, 28; 3:1 for 9; 8:1 for 12; 2:1:1 for 17, 19, 20; 3:1:1
for 16, 18; 3:2:1 for 21, 22, 25, 26; 3:3:1:1 for 23.

\newpage
\section{The table of other partitions}
 
{\bf Table 4. Embedding of some other partitions.} 

\[\begin{array}{||l||l||l|l||l|l||}\hline

29& De(L_5) & \alpha_3, Pyr_4 & ElDo& \frac{1}{2}{\bf Z}_4 & {\bf Z}_5 \\
30& De(D\mbox{-complex}) & \alpha_3, \sim \beta_3 & \mbox{triakis tr}
\alpha_3 & \frac{1}{2}{\bf Z}_5 & - \\
31& De(\mbox{Kelvin}) & \alpha_3, \beta_3 & RoDo, tw RoDo & - & - \\ 
32 & De(\mbox{Gr\"unbaum})& Prism_3 &Prism_6, BDS^* & \frac{1}{2}{\bf
Z}_5 & - \\
33 &De(\mbox{elong. Kelvin}) & \alpha_3, \beta_3, Prism_3 & RoDo-v & -
& - \\
34 & De(\mbox{elong. Gr\"unbaum})& Prism_3, \gamma_3 & \sim Prism_5 &
\frac{1}{2}{\bf Z}_5 & - \\
35 &P(S_1) & S_1 &  & {\bf Z}_3 & -  \\
36 &P(S_2) & S_2 &  & {\bf Z}_4 & -  \\
37 &P(S_3) & S_3 &  & {\bf Z}_5 & -  \\
38 & A\mbox{-19} & Prism_{\infty}& & {\bf Z}_2 & \\
39 &  & Prism_{\infty}& & {\bf Z}_2 & \\
40 & A \mbox{-20, $n$ even} & C_n \times P_{\bf Z} & & {\bf Z}_{\infty} & \\
   & A \mbox{-20, $n$ odd} & C_n \times P_{\bf Z} & & \frac{1}{2}
{\bf Z}_{\infty} & \\
41 & A\mbox{-22} & Aprism_{\infty} & & \frac{1}{2}{\bf Z}_3 & \\
42 & A\mbox{-23} & Prism_{\infty},Aprism_{\infty} & & \frac{1}{2}
{\bf Z}_3 & \\ 
43 & \parallel-type & \gamma_{3}, C_4 \times P_{\bf Z} & & {\bf Z}_3 &
\\
44 & \perp-type & \gamma_{3}, C_4 \times P_{\bf Z} & & {\bf Z}_3 &
\\
45 & chess-type & \gamma_{3}, C_4 \times P_{\bf Z} & & {\bf Z}_3 & \\ 
46 & A\mbox{-13}' & \alpha_3, \mbox{tr}\alpha_3 & R, \mbox{ twisted }R & - &-
\\ \hline
\end{array} \]
In Table 4 we group some other relevant partitions.
 Here $L_5$ denotes
a representative of the 5-th Fedorov's type (i.e. by the Voronoi polyhedron)  of
lattice in 3-space and $ElDo$ denotes its Voronoi polyhedron, called
elongated dodecahedron. Remaining four lattices appeared in Table 3 as
No 1 =$De({\bf Z}_3)=Vo({\bf Z}_3)$, No 5 =$De(A_3)$, No 2 =$Vo(A^*_3)$, 
No 3 =$De(A_2 \times {\bf Z}_1)$, No 4 =$Vo(A_2 \times {\bf Z}_1)$. 
Remark that $De(L_5)$ and $De(A_2 \times {\bf Z}_1)$ coincide as graphs,
but differ as partitions. 

In the notation $De$(Kelvin) below we consider any Kelvin packing by
$\alpha_3$ and 
$\beta_3$ (in proportion 2:1) which is {\em proper}, i.e. different from the
lattice $A_3$=fcc (face-centered lattice) and the bi-lattice hcp
(hexagonal closed packing). 
Any proper Kelvin partition, as well as partition $13'$ in \cite{An} 
(given as 46 in Table 4 and which Andreini
wrongly gave as uniform one), have exactly two vertex figures.
(The Voronoi tiles of tiling 46 are two rhombohedra: a rhombohedron,
say, $R$, i.e.
the cube contracted along a diagonal, and twisted $R$; both are
equivalent to $\gamma_3$.)  
The same is true for Gr\"unbaum partitions; see Section 5 below for those
notions and items 32--34 of Table 4. 

See Section 4 below for items 38--45 of Table 4. 
$D$-complex is the diamond bi-lattice; triakis tr$\alpha_3$ denotes
truncated $\alpha_3$ with $Pyr_3$ on each its triangular faces. 
Partitions 29 and 30 from Table 4 are both vertex-transitive, but they have
some non-Archimedean tiles: $Pyr_4$ for 29 and non-regular octahedron for 30.

The partitions 35, 36, 37 of Table 4 are all 3 non-normalizable
tilings of 3-space by convex parallelohedron, which where found in
\cite{Sh}. The polyhedra denoted by $S_1$, $S_2$ and $S_3$ are
centrally symmetric 10-hedra obtained by a decoration of the
paralelipiped. $S_1$ is equivalent to $\beta_3$ truncated on two
opposite vertices. $P^*(S_i)$ for $i=1,2,3$ are different partitions
of 3-space by non-convex bodies, but they have the same skeleton,
which is not 5-gonal. 

\section{Non-compact uniform partitions}
The non-compact uniform
partitions, introduced in subsections 19, 20, 22, 23 of \cite{An}, will be
denoted here A-19, A-20, A-22, A-23, respectively, and put in Table 4
as No's 38, 40, 41, 42. Denote by $Prism_\infty$ ($Antiprism_\infty$)
and $C_n \times P_{\bf Z}$ the $\infty$-sided prisms (antiprisms, 
respectively) and the cylinder on $C_n$. 

A-19 is obtained by putting $Prism_\infty$ on
($4^4$) and so its skeleton is ${\bf Z}_2$. We add, as item 39, the 
partition which differs from 38 only by another disposition of infinite
prisms {\em under} net ($4^4$), i.e. perpendicular to those above it. 

A-20 is obtained by putting the cylinders on ($4^4$), A-22 by
putting $Antiprism_\infty$ on ($3^6$) and A-23 by putting
$Prism_\infty$ and $Antiprism_\infty$ on ($3^3$.$4^2$).

 In subsection $20'$  \cite{An} mentions
also the partition into two half-spaces separated by some of 10
Archimedean (and one degenerated) nets, i.e. uniform plane
partitions. We can also take two parallel nets (say, $T$) and fill
the space between them by usual prisms (so, the skeleton will be 
direct product of the graph of $T$ and $K_2$) or, for $T = (4^4)$ or 
$(3^3.4^2)$, by a combination of usual and infinite prisms. Similar
uniform partitions are obtained if we will take an infinite number of
parallel nets $T$.

 The partitions 43, 44, 45 of Table 4 differ only by the disposition
of cubes and cylinders. In 43 the layers of cylinders stay
parallel ($\parallel$-type); in 44 they are perpendicular to the
cylinders of each previous layer. In 45 we see ($4^4$) as an infinite 
chess-board; cylinders stay on "white" squares while piles of cubes
stay on the "black" ones. 

By a decoration of $Prism_\infty$ in 38, 39,
one can get other non-compact uniform partitions.  

\section{Almost-uniform partitions} 
Call a normal partition of the 3-space into Platonic and Archimedean
polyhedra, {\em almost-uniform} if the group of symmetry is not
vertex-transitive but all vertex figures are congruent. 
Gr\"unbaum \cite{Gr} gave two infinite classes of such partitions and
indicated that he do not know 
other examples. In our terms, they called {\em elongated
proper Kelvin} and {\em elongated proper Gr\"unbaum} partitions. 
Kelvin and Gr\"unbaum partitions are defined uniquely by an infinite
binary sequence characterizing the way how layers follow each other. 
In Kelvin partition, the layers of $\alpha_3$ and $\beta_3$  follow 
each other in two
different ways (say, $a$ and $b$) while in Gr\"unbaum partition the layers
of $Prism_3$ follow each other in parallel or perpendicular mutual
disposition of heights.
Unproper Kelvin partitions give uniform partitions 5 and 24 for sequences
$...aaa...$ (or $...bbb...$) and $...ababab...$, respectively. 
Proper Kelvin and Gr\"unbaum partitions are not almost-uniform; 
there are even $\infty$-uniform
ones (take a non-periodic sequence).

Consider now elongations of those partitions, i.e. we add alternatively
the layers of $Prism_3$ for Kelvin and of cubes for Gr\"unbaum
partitions.

Remark that $RoDo-v$ (the Voronoi tile of partitions 25, 26, 33) can be
seen as a half of $RoDo$ cut in two, and that tw$RoDo$ is obtained from
$RoDo$ by a twist (a turn by $90^o$) of two halves.The Voronoi tiles for
proper Kelvin partition 31 are both $RoDo$ and tw$RoDo$ while only one of
them remains for two unproper cases 5, 24. Similarly, the Voronoi tile
of 34 (a special 5-prism) can be seen as a half of $Prism_6$ cut in two,
and $BDS^*$ can be seen as twisted $Prism_6$ in similar way. The Voronoi 
tiles for proper Gr\"unbaum partition 32 are both $Prism_6$ and
$BDS^*$ while only one of them remains for two unproper cases 3, 27.
 
Besides of two unproper cases 25, 26 (elongation of uniform 5, 24)
which are uniform, we have a continuum  of proper elongated Kelvin
partitions (denoted 33 in Table 4) which are almost-uniform. Among 
them there is a countable number of periodic partitions corresponding
to periodic $(a,b)$-sequences. Remaining continuum consists of aperiodic
tilings of 3-space by $\alpha_3$, $\beta_3$, $Prism_3$ with very simple
rule: each has unique Delaunay star consisting of 6 $Prism_3$ (put
together in order to form a 6-prism), 3 $\alpha_3$ and 3 $\beta_3$ (put
alternatively on 6 triangles subdividing the hexagon) and one $\alpha_3$
filling remaining space in the star. Each $b$ in the $(a,b)$-sequence,
defining such tiling, corresponds to the twist interchanging 3
$\alpha_3$ and 3 $\beta_3$ above (i.e. to the turn of the
configuration of 4 $\alpha_3$, 3 $\beta_3$ by $60^o$).

 Similar situation occurs for elongated Gr\"unbaum partitions. Besides
two uniform unproper cases 13, 28 (elongation of uniform 3, 27),
we have a continuum of proper elongated Gr\"unbaum partitions (denoted
34 in Table 4) which are almost-uniform. The aperiodic
$(a,b)$-sequences give a continuum of aperiodic tilings by $Prism_3$,
$\gamma_3$ with similar simple rule: unique Delaunay star consisting of
4 $\gamma_3$ (put together in order to form a 4-prism), 4 $Prism_3$
put on them and 2 $Prism_3$ filling remaining space in the star. Each
$b$ in $(a,b)$-sequence, defining the tiling, corresponds to a turn of
all configuration of 6 $Prism_3$ by $90^o$.

\section{Archimedean 4-polytopes}

Finite relatives of uniform partitions of 3-space are 4-dimensional
Archimedean polytopes, i.e. those having vertex-transitive group of
symmetry and whose cells are Platonic or Archimedean polyhedra and
prisms or antiprisms with regular faces. \cite{Con} enumerated all of
them:

1) 44 polytopes (others than prism on $\gamma_3$) obtained by Wythoff's
  kaleidoscope construction from 4-dimensional irreducible reflection
  (point) groups;

2) 17 prisms on Platonic (other than $\gamma_3$) and Archimedean solids;

3) Prism on $Antiprism_n$ with $n>3$;

4) A doubly infinity of polytopes which are direct products of two
  regular polygons (if one of polygons is a square, then we get a
  prisms on 3-dimensional prisms);

5) Gosset's semi-regular polytope called {\em snub 24-cell};

6) A new polytope, called {\em Grand Antiprism}, having 100 vertices 
(all from 600-cell), 300 cells $\alpha_3$ and 20 cells $Antiprism_5$ 
(those antiprisms form two interlocking tubes).

Using the fact that the direct product of two graphs is
$l_1$-embeddable if and only if each of them is, and the
characterization of embeddable Archimedean 
polyhedra in \cite{DS1} (see Table 1 above), we can decide on embeddability in cases
2)--4). In fact, the answer is "yes" always in cases 2)--4), except
prisms on tr$\alpha_3$, tr$\gamma_3$, Cuboctahedron, truncated
Icosahedron, truncated Dodecahedron and Icosidodecahedron, which all are
not 5-gonal. 

Now, the snub 24-cell embeds into $\frac{1}{2}H_{12}$ and the Grand 
Antiprism (as well as 600-cell itself) violates 7-gonal inequality, 
which is also necessary for embedding (see \cite{De}, \cite{DS1}).

\end{document}